\documentclass{article}
\newtheorem{Theorem}{Theorem}

\newtheorem{Corollary}[Theorem]{Corollary}

\newtheorem{Definition}[Theorem]{Definition}

\newtheorem{Axiom}{Axiom}

\begin{document}
\title{An axiomatic approach to diversity}
\author{Chris Dowden}
\maketitle
\setlength{\unitlength}{1cm}

\begin{abstract}
The topic of diversity is an interesting subject,
both as a purely mathematical concept and also for its applications to important real-life situations.
Unfortunately, although the meaning of diversity seems intuitively clear,
no precise mathematical definition exists.
In this paper, we adopt an axiomatic approach to the problem,
and attempt to produce a satisfactory measure.
\end{abstract}

\section{Introduction}

Over the last twenty years,
an important problem in conservation biology has been how best to measure the `diversity' of a set of organisms.
This is because diversity has emerged as a leading criterion when prioritising species to be saved from extinction.
The topic also has applications in a wide number of other fields, such as linguistics and economics,
but in this paper we examine it as a mathematical concept.

There are two distinct challenges.
The first is how to accurately evaluate the diversity of any two elements (e.g.~two organisms),
and the second is how to then use these pairwise-diversities, or `distances', to produce scores for sets of size greater than two.
It is the latter problem that we address here.

Biologists and economists have produced numerous papers 
(see~\cite{ber}--\cite{wei98} and the references therein)
investigating various different measures
that give reasonable approximations to diversity.
Some of these are very simple `rule of thumb' methods
(e.g.~minimum distance~\cite{bor}, maximum distance~\cite{ber}, average distance~\cite{war}, total distance~\cite{ber}),
while others are more elaborate
(e.g.~p-median~\cite{fai94}, phylogenetic diversity~\cite{fai92}, and Weitzman's diversity measure~\cite{wei92}, 
the latter two of which we shall shortly discuss).
Unfortunately, each of these is known to be imperfect,
in that they sometimes rank sets in a counter-intuitive order.

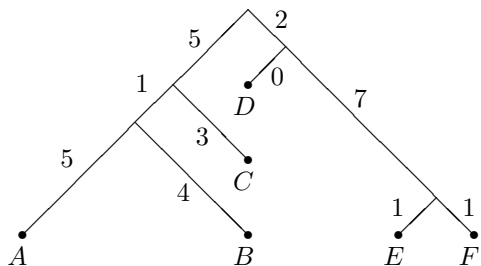
\begin{figure} [ht]
\setlength{\unitlength}{1cm}
\begin{picture}(12,4)(-4,-1.5)
\put(-1,-1){\line(1,1){3}}
\put(2,2){\line(1,-1){3}}
\put(1,1){\line(1,-1){1}}
\put(0.5,0.5){\line(1,-1){1.5}}
\put(2.5,1.5){\line(-1,-1){0.5}}
\put(4.5,-0.5){\line(-1,-1){0.5}}
\put(-1.2,-1.4){$A$}
\put(1.8,-1.4){$B$}
\put(1.8,-0.4){$C$}
\put(1.8,0.6){$D$}
\put(3.8,-1.4){$E$}
\put(4.8,-1.4){$F$}
\put(-0.5,-0.1){$5$}
\put(0.5,0.9){$1$}
\put(1.3,0.2){$3$}
\put(1.05,-0.55){$4$}
\put(1.2,1.5){$5$}
\put(2.35,1.75){$2$}
\put(2.3,1){$0$}
\put(3.4,0.7){$7$}
\put(3.9,-0.75){$1$}
\put(4.85,-0.75){$1$}
\put(-1,-1){\circle*{0.1}}
\put(2,-1){\circle*{0.1}}
\put(2,0){\circle*{0.1}}
\put(2,1){\circle*{0.1}}
\put(4,-1){\circle*{0.1}}
\put(5,-1){\circle*{0.1}}
\end{picture}
\caption{An example of phylogenetic diversity.} \label{PDtree}
\end{figure}

One of the most popular methods is `phylogenetic diversity' (\cite{fai92}).
Given the tree-like structure of evolutionary relationships,
phylogenetic diversity was developed for the specialised case when the pairwise-diversities induce a tree-metric,
with the score of a set of organisms being defined to be the length of the minimal subtree connecting them.
For example,
given the tree shown in Figure~\ref{PDtree},
the sets $\{A,C,D\}$ and $\{A,C,E\}$ score $16$ and $24$, respectively,
and so the latter would be considered as the more diverse.

Phylogenetic diversity has been observed to operate successfully in many scenarios.
However,
one problem that has been noted is that, in practice,
the pairwise-diversities will often not induce a tree-metric.
Also,
even with a tree-metric,
sets can still sometimes appear to be ranked in a slightly undesirable order.
For example,
in Figure~\ref{PDtree} the set $\{A,E,F\}$ would score more than the set $\{A,D,F\}$ ($22$ compared to $21$),
even though $D$ is very different from both $A$ and $F$,
while $E$ is very similar to $F$.
Indeed, adding $D$ to the set $\{A,F\}$ does not increase the phylogenetic diversity score at all,
which seems counter-intuitive.

Another prominent diversity measure is that investigated by Weitzman in his influential paper `On Diversity' (\cite{wei92}).
Here, the diversity $V$ of a set $S$ is defined recursively by the formula
\begin{displaymath}
V(S) = \max_{i} \{V(S \setminus i) + d(i,S \setminus i)\},
\end{displaymath}
where $d(i,S \setminus i)$ denotes the minimum distance between $i$ and an element of $S \setminus i$.
For example,
in Figure~\ref{PDtree} we would obtain
$V(\{A,E,F\}) = \max \{21+2,21+2,2+21\} = 23$
and $V(\{A,D,F\}) = \max \{13+8,21+8, 8+13\} = 29$,
which seems reasonable.

Again,
this measure can be observed to operate successfully in many cases,
but unfortunately
there are still circumstances where the formula can produce scores that appear slightly imperfect.
For example,
it seems natural that the set $\{A,B,F\}$
should be considered slightly more diverse than the set $\{A,C,F\}$,
since both $B$ and $C$ are distance $9$ from $A$ 
but $B$ is further than  $C$ from $F$.
However,
it can be checked that Weitzman's recursion gives $V(\{A,B,F\})=30=V(\{A,C,F\})$,
i.e.~there is no preference between the two sets.

It is the object of this paper to introduce a new axiomatic approach to diversity,
in an attempt to produce a measure that never disagrees with intuition.
Furthermore, we shall only assume that the pairwise-diversities satisfy the properties of a metric,
and not necessarily a tree-metric.

Before we make a start,
it is important to point out that the diversity of the set of features (or genes)
contained by a set of organisms is often of interest when making conservation decisions
(see e.g.~\cite{neh02} or~\cite{wei98}),
rather than just the diversity of the set of organisms themselves.
Although clearly closely related,
these two criteria are not identical,
and it is not always the case that choosing the most diverse set of organisms equates
to choosing the organisms that provide the most diverse set of features.

For example,
if three organisms are represented by the sets 
$A=\{f_{1},f_{2},f_{3},f_{4},f_{5}\}$,
$B=\{f_{1},f_{2},f_{3},f_{4},f_{6}\}$
and $C=\{f_{1}\}$
(where the $f_{i}$'s denote various features/genes),
then the set of features contained by $A$ and $B$ (i.e.~$A \cup B$)
would be considered to have slightly greater diversity 
than the set of features contained by $A$ and $C$ (i.e.~$A \cup C$),
even though the set $\{A,B\}$ would probably be considered to be less diverse
than the set $\{A,C\}$.
As a more trivial example,
the set $A$ might score highly for diversity,
even though the set $\{A\}$ would have zero diversity.

Finally, it is worth remarking that,
even without the biological motivation,
the question under discussion in this paper seems very natural ---
given two collections of points in a metric space,
which is the more spread out?
It is perhaps surprising that the topic seems not to have been widely investigated by mathematicians.

The remainder of the first half of this paper is divided into two sections.
In the first, Section~\ref{axioms},
we propose four basic axioms that a diversity measure should satisfy;
in the second,
Section~\ref{basicdm}, 
we then present some measures that fulfill all these requirements
(the first to be produced).

In the second half of the paper,
we then discuss one further potential axiom.
In Section~\ref{equidistance},
we present a limited version,
and show how to modify our earlier measures to also satisfy this property;
in Section~\ref{strongsection},
we investigate a seemingly natural `strong' version,
which we shall see actually produces an intriguing contradiction;
and in Section~\ref{nice measure},
we outline a more careful formulation of the desired property,
and describe a pretty measure that appears to always give nice results.

\section{Axioms} \label{axioms}

Throughout the rest of this paper,
we shall assume that we are given a complete weighted graph,
where the edge-weights denote the (given) pairwise-diversities of the vertices
(this is slightly different from the tree-like structure of Figure~\ref{PDtree}).
For ease of expression,
we make the following definitions:

\begin{Definition}
Let us use $D(S)$ to denote the diversity of a set $S$,
and let us say
\emph{the distance between $x$ and $y$}
to mean $D(\{x,y\})$.
We will assume only that these distances
(i.e.~the pairwise-diversities)
satisfy the properties of a metric.
\end{Definition}

Our aim is to construct a way to use the distances to give a score for the overall diversity of any subset of our collection of vertices.
To that end, we will spend this section proposing four axioms.

We start with three properties that are hoped to be intuitively natural:

\begin{Axiom} \label{adding}
For any non-empty set of vertices $S$,
we have $D(S \cup \{x\}) \geq D(S)$ for all $x$, with equality if and only if $x \in S$.
\end{Axiom}

\begin{Axiom} \label{increasing}
For any two vertex-sets $S = \{s_{1},s_{2}, \ldots, s_{n}\}$ and $T = \{t_{1},t_{2}, \ldots, t_{n}\}$
satisfying $D(\{t_{i},t_{j}\}) \geq D(\{s_{i},s_{j}\})$ for all $i$ and $j$, we have $D(T) \geq D(S)$,
with equality if and only if $D(\{t_{i},t_{j}\}) = D(\{s_{i},s_{j}\})$ for all $i$ and $j$.
\end{Axiom}

\begin{Axiom} \label{continuity}
\emph{Continuity}.
Given any set of vertices $S = \{s_{1},s_{2}, \ldots, s_{n}\}$ and any $\epsilon > 0$,
there exists a $\delta (S, \epsilon) > 0$
such that, for any set of vertices $T = \{t_{1},t_{2}, \ldots, t_{n}\}$
satisfying $|D(\{t_{i},t_{j}\}) - D(\{s_{i},s_{j}\})| < \delta$ for all $i$ and $j$,
we have $|D(T)-D(S)| < \epsilon$.
\end{Axiom}

It is worth observing that two other desirable properties follow automatically from these axioms.
First, note that Axiom~\ref{adding} implies $D(\{x\}) = D(\{x,x\})$, and hence
(since we assume that the distances satisfy the properties of a metric):

\begin{Corollary}
$D(\{x\})=0$ for all $x$.
\end{Corollary}
Secondly, it follows from applying Axiom~\ref{continuity} to the sets
$S = \{s_{1},s_{2}, \ldots, s_{n},s_{n}\}$ and $T = \{s_{1},s_{2}, \ldots, s_{n},x\}$
(and using the triangle inequality) that we have continuity when adding a new vertex:

\begin{Corollary}
Given a set of vertices $S = \{s_{1},s_{2}, \ldots, s_{n}\}$ and an $\epsilon > 0$,
there exists a $\delta (S, \epsilon) > 0$
such that, for any vertex $x$ satisfying $D(\{x,s_{n}\}) < \delta$,
we have $D(S \cup \{x\}) < D(S) + \epsilon$.
\end{Corollary}

Our fourth axiom is motivated by the principle that consistent results ought to be obtained regardless of differences in the scale used
to measure the original distances.
For example, if we wish to compare the diversity of the locations of stars in two different galaxies,
then the resultant ranking should not depend on whether the distances were measured in light-years or kilometres.
In other words, multiplying all our original distances by some constant $c$ should not affect
whether or not $D(S) > D(T)$ for any sets $S$ and $T$:

\begin{Axiom} \label{scaling}
\emph{Scaling.}
Given four sets of vertices $S= \{s_{1},s_{2}, \ldots, s_{n}\}$,
$S^{\prime} = \{s_{1}^{\prime},s_{2}^{\prime}, \ldots, s_{n}^{\prime}\}$,
$T = \{t_{1},t_{2}, \ldots, t_{k}\}$
and $T^{\prime} = \{t_{1}^{\prime},t_{2}^{\prime}, \ldots, t_{k}^{\prime}\}$,
if $D(\{s_{i}^{\prime},s_{j}^{\prime}\}) = cD(\{s_{i},s_{j}\})$ for all $i$ and $j$
and $D(\{t_{i}^{\prime},t_{j}^{\prime}\}) = cD(\{t_{i},t_{j}\})$ for all $i$ and $j$, 
for some constant $c>0$,
then $D(S^{\prime}) > D(T^{\prime})$ if and only if $D(S) > D(T)$.
\end{Axiom}

By considering the case when $|T|=2$,
this implies the following:

\begin{Corollary} 
Given two sets $S = \{s_{1},s_{2}, \ldots, s_{n}\}$ and $S^{\prime} = \{s_{1}^{\prime},s_{2}^{\prime}, \ldots, s_{n}^{\prime}\}$,
if $D(\{s_{i}^{\prime},s_{j}^{\prime}\}) = cD(\{s_{i},s_{j}\})$ for all $i$ and $j$, for some constant $c>0$,
then $D(S^{\prime}) = cD(S)$.
\end{Corollary}

In the following section,
we shall present measures that can be be shown to satisfy all four of our basic axioms.
However,
there is then also an additional `equidistance' axiom
that we will discuss extensively in Sections~\ref{equidistance}--\ref{nice measure}.

\section{New diversity measures} \label{basicdm}

In the previous section, we proposed four axioms that a diversity measure should satisfy.
Although the problem seems fairly natural,
it is surprisingly difficult to construct a measure that fulfills all these requirements
(indeed, every method examined in the existing literature seems to fail either Axiom~\ref{adding} or Axiom~\ref{increasing}).
However, in this section we shall now present a simple system
for obtaining measures that do satisfy all four axioms.
 
The purpose of introducing these measures is twofold.
Firstly,
to show that it is indeed possible to satisfy all of Axioms~\ref{adding}--\ref{scaling};
secondly,
to show (in Section~\ref{equidistance})
that measures satisfying Axioms~\ref{adding}--\ref{scaling}
can still produce counter-intuitive results,
and that a fifth axiom is consequently necessary.

\begin{Definition}
Given a real-valued function $f$
defined on all vertex-sets of size at least three,
let us define the measure $D_{f}$ recursively (from our given distances) by using the equation
\begin{equation} \label{basic}
D_{f}(S) = f(S) + \max_{T \subset S, |T|=|S|-1} \{ D_{f} (T) \}
\end{equation}
for all vertex-sets $S$ of size greater than two.
Let us call the function $f$ \emph{suitable} if:
(a) $f$ is a continuous function of the distances;
(b) if any of the distances are $0$, then $f=0$;
(c) if none of the distances are $0$, then $f$ is strictly positive and is a monotonically strictly increasing function of the distances;
and (d) $f$ is `scale-invariant' in the sense of Axiom~\ref{scaling}.
\end{Definition}

As an example of a `suitable' function, 
we could choose $f(\{s_{1},s_{2}, \ldots, s_{n}\})$ to be
$\left( \prod_{1 \leq i < j \leq n} D(\{s_{i},s_{j}\}) \right)^{\frac{1}{\left(^{n}_{2}\right)}}$
or $\left( \sum_{1 \leq i < j \leq n} \frac{1}{D(\{s_{i},s_{j}\})} \right)^{-1}$,
or any linear combination of these.

We will now see that $D_{f}$ satisfies the axioms if $f$ is suitable:

\begin{Theorem} \label{satisfaction}
The diversity measure $D_{f}$ defined in equation~(\ref{basic}) satisfies Axioms~\ref{adding}--\ref{scaling} if the function $f$ is suitable.
\end{Theorem}
\textbf{Proof}
It is immediately clear by induction that $D_{f}$ satisfies Axioms~\ref{increasing}--\ref{scaling},
so it only remains to show that Axiom $1$ is satisfied.
To do this,
we need to prove that
(i) adding a vertex already in the set does not alter the score,
and (ii) adding a vertex not already in the set strictly increases the score.

We shall proceed by induction.
Suppose that (i) and (ii) both hold when adding a vertex to any set of size less than $k$
(note that the base step is a direct consequence of the fact that the distances satisfy the properties of a metric),
and let us now consider the case when we are adding a vertex $s_{k+1}$ to a set $S=\{s_{1},s_{2}, \ldots, s_{k}\}$ of size exactly $k$.

First, let us work towards (i) by supposing (without loss of generality) that $s_{k+1} = s_{k}$.
By part (b) of the definition of suitability,
we then have $f\big(S \cup~\!\!\{s_{k+1}\}\big) = 0$
and so $D_{f}\big(S \cup \{s_{k+1}\}\big)  = \max_{i \leq k+1} \left\{D_{f}\Big(\big(S \cup \left\{s_{k+1}\right\}\big) \setminus s_{i}\Big)\right\}$.
Hence, it suffices to prove that
$D_{f}\Big(\big(S \cup \{s_{k+1}\}\big) \setminus s_{i}\Big)$
is maximised by taking $i \in \{k,k+1\}$.
But note that, for $i<k$,
the induction hypothesis implies
$D_{f}\Big(\big(S \cup \{s_{k+1}\}\big) \setminus s_{i}\Big) = D_{f}(S \setminus s_{i}) \leq D_{f}(S)$,
and so we are done.

Now let us work towards (ii) by supposing that $s_{k+1} \notin S$.
If the vertices of $S$ are all distinct,
then $f(S \cup \{s_{k+1}\}) > 0$ and the result is clear.
If not, then let $S^{\prime}$ denote a maximally sized subset of $S$ with vertices that are all distinct.
By the induction hypothesis,
we have $D_{f}(S^{\prime} \cup \{s_{k+1}\}) > D_{f}(S^{\prime})$.
But note that, by a combination of the induction hypothesis and (i),
the left-hand side is
$D_{f}(S \cup \{s_{k+1}\})$
and the right-hand side is $D_{f}(S)$.~
$\phantom{qwerty}
\setlength{\unitlength}{.25cm}
\begin{picture}(1,1)
\put(0,0){\line(1,0){1}}
\put(0,0){\line(0,1){1}}
\put(1,1){\line(-1,0){1}}
\put(1,1){\line(0,-1){1}}
\end{picture}$ \\

One particular choice for a suitable function would be to take
\begin{equation} \label{nice}
f(\{s_{1},s_{2}, \ldots, s_{n}\}) = \frac{\left(^{n}_{2}\right)}{\sum_{1 \leq i < j \leq n} \frac{1}{D(\{s_{i},s_{j}\})}},
\end{equation}
and we shall refer back to this.
Note that this particular function has the aesthetically pleasing property
that it will always be equal to $1$ for the `regular' case
when the distances are all $1$,
and hence $D_{f}(S)$ will be equal to $|S|-1$ for this case.

\section{Three-vertex equidistance} \label{equidistance}

In the previous section,
we saw a scheme for generating diversity measures that satisfy Axioms~\ref{adding}--\ref{scaling}.
However, as briefly mentioned earlier, there is also a fifth axiom that is necessary --- that of \emph{equidistance}.
In this section, we shall explain why such an axiom is desirable,
and define it for the specific case when our graphs have exactly three vertices,
which we shall see only requires a small modification to our previous measures.
In Sections~\ref{strongsection} and \ref{nice measure},
we will then investigate how to extend the concept to graphs of arbitrary size.

Let us imagine that we have two vertices $x$ and $y$
that are distance $1$ apart,
and that we wish to add one more vertex to this set.
Suppose that we are free to choose any element from
$\{z: D(\{x,z\}) + D(\{y,z\}) = 2\}$.
It seems natural that the overall diversity score ought to be greater the more equidistant the new vertex is between $x$ and $y$.
Unfortunately,
this is actually not true for the diversity measure defined at the end of the last section,
where we use the suitable function of equation~(\ref{nice}) in recursion~(\ref{basic}),
since we know that the regular unit triangle scores $2$,
whereas the triangle with lengths $1$, $\frac{1}{2}$ and $\frac{3}{2}$ scores
$\frac{3}{1+2+\frac{2}{3}} + \max \left\{ 1, \frac{1}{2}, \frac{3}{2} \right\} = \frac{9}{11} + \frac{3}{2} > 2$.
This example establishes the need for a new axiom:

\begin{Axiom} \label{basic equidistance}
\emph{Three-vertex equidistance}.
Given two sets $S = \{s_{1}, s_{2}, s_{3} \}$ and $T = \{t_{1}, t_{2}, t_{3} \}$,
if $D(\{t_{1},t_{2}\}) = D(\{s_{1},s_{2}\})$
and $D(\{t_{1},t_{3}\}) + D(\{t_{2},t_{3}\}) = D(\{s_{1}, s_{3}\}) + D(\{s_{2}, s_{3}\}) = \lambda$,
but $\left|D(\{t_{1},t_{3}\}) - \frac{\lambda}{2}\right| < \left|D(\{s_{1},s_{3}\}) - \frac{\lambda}{2}\right|$
(and  
$\left|D(t_{2},t_{3}) - \frac{\lambda}{2}\right| <
\left|D(\{s_{2},s_{3}\}) - \frac{\lambda}{2}\right|$),
then $D(T) > D(S)$.
\end{Axiom}

One way to approach the problem of trying to satisfy Axiom~\ref{basic equidistance}
would be to find a suitable function $f$ for which the partial derivative with respect to the maximum distance
(when the total distance is fixed)
is always less than $-1$,
thus offsetting the contribution to $D_{f}$ of the $\max_{T \subset S : |T|=2} \{D_{f}(T)\}$ term.

However, a neater solution is to instead develop a separate diversity measure for sets of size three that does satisfy Axiom~\ref{basic equidistance}
(and also Axioms~\ref{adding}--\ref{scaling})
and then simply use the recursion of equation~(\ref{basic}) on this.
This prompts the following definition:

\begin{Definition}
Given a real-valued function $g(\{s_{1},s_{2},s_{3}\})$
on sets of size three,
and a real-valued function $f$ on sets of size greater than three,
let us define $D_{f,g}(\{s_{1},s_{2},s_{3}\})$ by the equation
\begin{equation} \label{Dfg1}
D_{f,g}(\{s_{1},s_{2},s_{3}\}) = g(\{s_{1},s_{2},s_{3}\}) + \frac{1}{2} \sum_{1 \leq i<j \leq 3} D(\{s_{i},s_{j}\}),
\end{equation}
and let us then define $D_{f,g}(S)$ recursively by using equation~(\ref{basic}), i.e.
\begin{equation} \label{Dfg2}
D_{f,g}(S) = f(S) + \max_{T \subset S, |T|=|S|-1} \{ D_{f,g} (T) \}
\end{equation}
for all vertex-sets $S$ of size greater than three.

Let us say that the function $g$ \emph{satisfies Axiom~\ref{basic equidistance}}
if given any two sets 
$S = \{s_{1}, s_{2}, s_{3} \}$ and $T = \{t_{1}, t_{2}, t_{3} \}$
for which $D(\{t_{1},t_{2}\}) = D(\{s_{1},s_{2}\})$
and $D(\{t_{1},t_{3}\}) + D(\{t_{2},t_{3}\}) = D(\{s_{1}, s_{3}\}) + D(\{s_{2}, s_{3}\}) = \lambda$,
but $\left|D(\{t_{1},t_{3}\}) - \frac{\lambda}{2}\right| < \left|D(\{s_{1},s_{3}\}) - \frac{\lambda}{2}\right|$
(and  
$\left|D(t_{2},t_{3}) - \frac{\lambda}{2}\right| <
\left|D(\{s_{2},s_{3}\}) - \frac{\lambda}{2}\right|$),
we always have $g(T) > g(S)$.
\end{Definition}

\begin{Theorem} \label{3vertexsectionthm}
The diversity measure $D_{f,g}$ defined in equations~(\ref{Dfg1}) and~(\ref{Dfg2})
satisfies Axioms~\ref{adding}--\ref{basic equidistance}
if the functions $f$ and $g$ are suitable
and the function $g$ satisfies Axiom~\ref{basic equidistance}.
\end{Theorem}
\textbf{Proof}
If $g$ is a suitable function that also satisfies Axiom~\ref{basic equidistance},
then it is simple to see that $D_{f,g}(\{s_{1},s_{2},s_{3}\})$ satisfies Axioms~\ref{increasing}--\ref{basic equidistance},
while Axiom~\ref{adding} (for the case when we are comparing sets of size three with sets of size two)
follows from using the triangle inequality on the term 
$ \frac{1}{2} \sum_{1 \leq i<j \leq 3} D(\{s_{i},s_{j}\})$.

It can then be checked that the proof of Theorem~\ref{satisfaction} will still hold,
and so $D_{f,g}(S)$ will also satisfy the axioms for sets of size greater than three.~
$\phantom{qwerty}
\setlength{\unitlength}{.25cm}
\begin{picture}(1,1)
\put(0,0){\line(1,0){1}}
\put(0,0){\line(0,1){1}}
\put(1,1){\line(-1,0){1}}
\put(1,1){\line(0,-1){1}}
\end{picture}$ \\

If we use the function $f$ defined in equation~(\ref{nice}),
then an aesthetically pleasing choice for the function $g$
(if we again wish to ensure that
$D_{f,g}(S)=|S|-1$
for the `regular' case when the distances are all $1$)
is to take $g(\{s_{1},s_{2},s_{3}\}) = \frac{3}{2} \left( \sum_{1 \leq i < j \leq 3} \frac{1}{D(\{s_{i},s_{j}\})} \right)^{-1}$,
so that 
\begin{eqnarray}
D_{f,g}(\{s_{1},s_{2},s_{3}\}) = \frac{3}{2} \left( \sum_{1 \leq i < j \leq 3} \frac{1}{D(\{s_{i},s_{j}\})} \right)^{-1}
+ \frac{1}{2} \sum_{1 \leq i<j \leq 3} D(\{s_{i},s_{j}\}). \label{size3formula}
\end{eqnarray}
We shall refer back to this later,
after Theorem~\ref{3setthm}.

\section{`Strong' equidistance} \label{strongsection}

In the previous section,
we defined equidistance only for the case of graphs with three vertices.
In this section,
we shall now present a seemingly natural way to extend this idea to larger graphs,
before showing that the resultant axiom would actually not be compatible with Axioms~\ref{adding}--\ref{continuity}!
In Section~\ref{nice measure},
we will then investigate matters further.

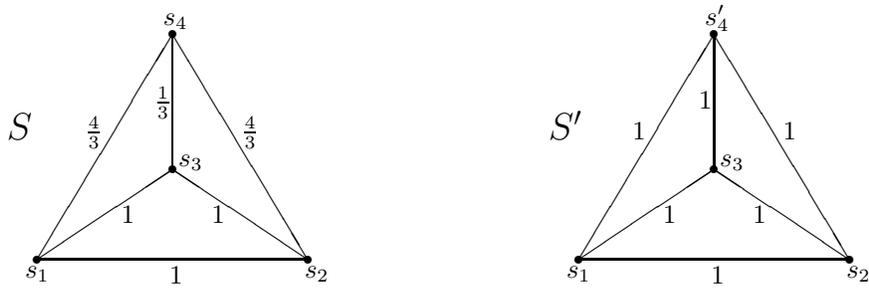
\begin{figure} [ht]
\setlength{\unitlength}{0.8cm}
\begin{picture}(10,4.5)(-0.75,-0.3)
\put(0,0){\line(1,0){4.5}}
\put(0,0){\line(3,5){2.25}}
\put(2.25,3.75){\line(3,-5){2.25}}
\put(2.25,3.75){\line(0,-1){2.25}}
\put(0,0){\line(3,2){2.25}}
\put(4.5,0){\line(-3,2){2.25}}
\put(-0.2,-0.3){$s_{1}$}
\put(4.45,-0.3){$s_{2}$}
\put(2.1,3.9){$s_{4}$}
\put(2.35,1.55){$s_{3}$}
\put(2.2,-0.4){$1$}
\put(2.9,0.6){$1$}
\put(1.4,0.6){$1$}
\put(0.8,2){$\frac{4}{3}$}
\put(3.4,2){$\frac{4}{3}$}
\put(1.95,2.5){$\frac{1}{3}$}
\put(9,0){\line(1,0){4.5}}
\put(9,0){\line(3,5){2.25}}
\put(11.25,3.75){\line(3,-5){2.25}}
\put(11.25,3.75){\line(0,-1){2.25}}
\put(9,0){\line(3,2){2.25}}
\put(13.5,0){\line(-3,2){2.25}}
\put(8.8,-0.3){$s_{1}$}
\put(13.45,-0.3){$s_{2}$}
\put(11.1,3.9){$s_{4}^{\prime}$}
\put(11.35,1.55){$s_{3}$}
\put(11.2,-0.4){$1$}
\put(11.9,0.6){$1$}
\put(10.4,0.6){$1$}
\put(9.9,2){$1$}
\put(12.4,2){$1$}
\put(11,2.5){$1$}
\put(-0.5,2){\Large{$S$}}
\put(8.5,2){\Large{$S^{\prime}$}}
\put(0,0){\circle*{0.125}}
\put(4.5,0){\circle*{0.125}}
\put(2.25,1.5){\circle*{0.125}}
\put(2.25,3.75){\circle*{0.125}}
\put(9,0){\circle*{0.125}}
\put(13.5,0){\circle*{0.125}}
\put(11.25,1.5){\circle*{0.125}}
\put(11.25,3.75){\circle*{0.125}}
\end{picture}
\caption{Two four-vertex sets, $S$ and $S^{\prime}$.} \label{c1}
\end{figure}

To see that something more general 
than just the three-vertex rule of Axiom~\ref{basic equidistance} is needed,
consider the two four-vertex sets, $S$ and $S^{\prime}$, depicted in Figure~\ref{c1}.
It seems intuitive that $S^{\prime}$ should be considered as more diverse than $S$,
as the position of $s_{4}^{\prime}$ is equidistant in relation to $s_{1}$, $s_{2}$ and $s_{3}$.
Unfortunately, the diversity measure that we defined in equation~(\ref{size3formula}) gives $D_{f,g}(S^{\prime}) = 3$
and $D_{f,g}(S) = \frac{6}{\frac{3}{4} + \frac{3}{4} + 3 + 1 + 1 +1} + D_{f,g}(\{s_{1},s_{2},s_{4}\})
= \frac{4}{5} + \frac{3}{2} \frac{1}{\frac{3}{4} + \frac{3}{4} + 1} + \frac{1}{2} \left( \frac{4}{3} + \frac{4}{3} + 1 \right) = \frac{97}{30} > 3$.

One natural way to extend Axiom~\ref{basic equidistance} to any number of vertices seems to be the following: \\
\\
\textbf{`Axiom' $\mathbf{5^{\prime}}$}
Strong equidistance.
\emph{Given two sets $S = \{s_{1}, s_{2}, \ldots, s_{n} \}$ and $T = \{t_{1}, t_{2}, \ldots, t_{n} \}$,
if $D(\{t_{i},t_{j}\}) = D(\{s_{i},s_{j}\})$ for all $i,j < n$
and $\sum_{i<n}D(\{t_{i},t_{n}\}) = \sum_{i<n} D(\{s_{i}, s_{n}\}) = \lambda$,
but 
$\left|D(\{t_{i},t_{n}\}) - \frac{\lambda}{n-1}\right| \stackrel{(*)}{\leq} \left|D(\{s_{i},s_{n}\}) - \frac{\lambda}{n-1}\right|$
for all $i<n$,
then $D(T) \geq D(S)$,
with equality if and only if there is equality in (*) for all $i<n$.} \\

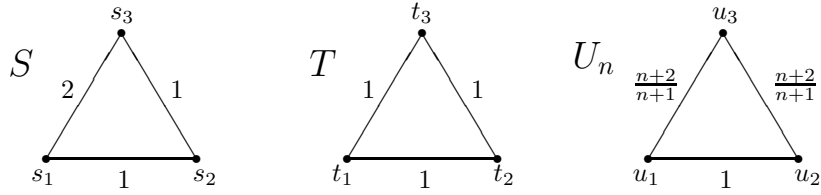
\begin{figure} [ht]
\setlength{\unitlength}{1cm}
\begin{picture}(10,2.55)(-1,-0.4)
\put(0,0){\line(1,0){2}}
\put(0,0){\line(3,5){1}}
\put(2,0){\line(-3,5){1}}
\put(4,0){\line(1,0){2}}
\put(4,0){\line(3,5){1}}
\put(6,0){\line(-3,5){1}}
\put(8,0){\line(1,0){2}}
\put(8,0){\line(3,5){1}}
\put(10,0){\line(-3,5){1}}
\put(0,0){\circle*{0.1}}
\put(2,0){\circle*{0.1}}
\put(1,1.666667){\circle*{0.1}}
\put(-0.2,-0.3){$s_{1}$}
\put(1.95,-0.3){$s_{2}$}
\put(0.85,1.85){$s_{3}$}
\put(0.2,0.8){$2$}
\put(0.95,-0.4){$1$}
\put(1.65,0.8){$1$}
\put(-0.5,1.1){\Large{$S$}}
\put(4,0){\circle*{0.1}}
\put(6,0){\circle*{0.1}}
\put(5,1.666667){\circle*{0.1}}
\put(3.8,-0.3){$t_{1}$}
\put(5.95,-0.3){$t_{2}$}
\put(4.85,1.85){$t_{3}$}
\put(4.2,0.8){$1$}
\put(4.95,-0.4){$1$}
\put(5.65,0.8){$1$}
\put(3.5,1.1){\Large{$T$}}
\put(8,0){\circle*{0.1}}
\put(10,0){\circle*{0.1}}
\put(9,1.666667){\circle*{0.1}}
\put(7.8,-0.3){$u_{1}$}
\put(9.95,-0.3){$u_{2}$}
\put(8.85,1.85){$u_{3}$}
\put(7.8,0.9){$\frac{n+2}{n+1}$}
\put(8.95,-0.4){$1$}
\put(9.65,0.9){$\frac{n+2}{n+1}$}
\put(7,1.2){\Large{$U_{n}$}}
\end{picture}
\caption{The sets $S$, $T$ and $U_{n}$ in the proof of Theorem~\ref{contradiction}.} \label{c2}
\end{figure}

Unfortunately,
as we shall now see,
it turns out that this strong version actually leads to inconsistencies with our earlier axioms!:

\begin{Theorem} \label{contradiction}
`Axiom'~$5^{\prime}$ is inconsistent with Axioms~\ref{adding}--\ref{continuity}.
\end{Theorem}
\textbf{Proof}
Let the vertex-sets $S= \{s_{1},s_{2},s_{3}\}$, $T=\{t_{1},t_{2},t_{3}\}$ and $U_{n}=\{u_{1},u_{2},u_{3}\}$ be as shown in Figure~\ref{c2}.
By Axiom~\ref{increasing},
we have $D(S) > D(T)$ and so, by continuity (Axiom~\ref{continuity}),
there exists a $k$ such that $D(S) > D(U_{k})$.

Now consider the set $U_{k}^{\prime} = \{ u_{1},u_{2}, \ldots, u_{k+2} \} \supset U_{k}$
constructed from $U_{k}$ by setting $D(\{u_{2},u_{l}\}) = 0$ for all $l \geq 4$
(i.e.~adding in $k-1$ extra copies  of $u_{2}$)
and, similarly, the set $S^{\prime} = \{s_{1},s_{2}, \ldots, s_{k+2}\}$ constructed from $S$
by setting $D(\{s_{2},s_{l}\}) = 0$ for all $l \geq 4$
(i.e.~adding in $k-1$ extra copies of $s_{2}$).

If we assume the strong equidistance of `Axiom'~$5^{\prime}$,
then $D(U_{k}^{\prime}) > D(S^{\prime})$.
But, by Axiom~\ref{adding}, $D(U_{k}^{\prime}) = D(U_{k})$ and $D(S^{\prime}) = D(S)$.
Hence, we find that $D(U_{k}) > D(S)$,
and so we have a contradiction.~
$\phantom{qwerty}
\setlength{\unitlength}{.25cm}
\begin{picture}(1,1)
\put(0,0){\line(1,0){1}}
\put(0,0){\line(0,1){1}}
\put(1,1){\line(-1,0){1}}
\put(1,1){\line(0,-1){1}}
\end{picture}$ \\

Note that
(with a bit of care to ensure that the triangle inequality is not violated during the proof)
a form of this example still produces a contradiction even if we alter the strong equidistance `axiom' to include the condition
$D(\{t_{i},t_{n}\})=D(\{s_{i},s_{n})\}$ for all $i \geq 3$ as well as
$\sum_{i<n}D(\{t_{i},t_{n}\}) = \sum_{i<n}D(\{s_{i},s_{n}\})$.
This seems very surprising!

\section{Symmetry considerations and a proposed diversity measure} \label{nice measure}

In this section,
we shall conclude our exploration of diversity axioms with a brief mention of symmetry considerations,
and we shall present some speculative work concerning a final diversity measure that appears satisfactory.

Let us recall our four-vertex example of Figure~\ref{c1}.
The intuition that it was desirable for the fourth vertex to be equidistant in relation to $s_{1}$, $s_{2}$ and $s_{3}$
perhaps stems from the fact that these other three vertices were all in symmetric positions to begin with
(note that Axiom~\ref{basic equidistance} is also an application of this intuition,
since every set of size two is automatically symmetric).
Hence, it seems sensible that a truly satisfactory diversity measure should have to take into account such symmetry considerations
(a lengthy discussion of such matters,
including proposed `symmetric equidistance' axioms,
is given in an earlier version of this paper~(\cite{dow})).

Unfortunately,
the measures given in Section~\ref{equidistance} seem irreparably distorted by the $\max \{D_{f,g}(T)\}$ term,
which was used to satisfy Axiom~\ref{adding}.
However, we shall now present a possible alternative that appears to give nice results without employing such an expression.

\begin{Definition}
Given a set $S = \{s_{1},s_{2}, \ldots, s_{n}\}$,
let $p_{kl} = \frac{ \frac{1}{D(\{s_{k},s_{l}\})}}{\sum_{1\leq i<j \leq n} \frac{1}{D(\{s_{i},s_{j}\})}}$ for all $k \neq l$,
and define $D$ recursively by the equation
\begin{eqnarray}
D(S) = \sum_{1 \leq k<l \leq n} p_{kl} \Big(D\left(\left\{s_{k},s_{l}\right\}\right) + D(S_{kl})\Big), \label{nice diversity}
\end{eqnarray}
where $S_{kl}$ denotes the set formed from $S$ by `merging' $s_{k}$ and $s_{l}$ into a new vertex $s_{kl}$
and setting $D(\{s_{kl},s_{i}\}) = \frac{D(\{s_{k},s_{i}\}) + D(\{s_{l},s_{i}\})}{2}$ for all $i$
(it can be checked that the distances in $S_{kl}$ will still satisfy the properties of a metric).
\end{Definition}

\begin{figure} [ht]
\setlength{\unitlength}{1cm}
\begin{picture}(10,2.4)(-0.6,-0.3)
\put(0,0){\line(1,0){2}}
\put(0,0){\line(3,5){1}}
\put(2,0){\line(-3,5){1}}
\put(0,0){\circle*{0.1}}
\put(2,0){\circle*{0.1}}
\put(1,1.666667){\circle*{0.1}}
\put(-0.2,-0.3){$s_{1}$}
\put(1.95,-0.3){$s_{2}$}
\put(0.85,1.85){$s_{3}$}
\put(0.2,0.8){$2$}
\put(0.95,-0.4){$4$}
\put(1.65,0.8){$3$}
\put(-0.5,1.1){\Large{$S$}}
\put(4.5,0){\line(0,1){1.666667}}
\put(4.5,0){\circle*{0.1}}
\put(4.5,1.666667){\circle*{0.1}}
\put(4.35,-0.3){$s_{12}$}
\put(4.35,1.85){$s_{3}$}
\put(4.6,0.8){$2.5$}
\put(3.5,1.1){\Large{$S_{12}$}}
\put(8,0){\line(-2,1){1.5}}
\put(6.5,0.75){\circle*{0.1}}
\put(8,0){\circle*{0.1}}
\put(6.3,0.4){$s_{13}$}
\put(7.95,-0.3){$s_{2}$}
\put(7.2,0.5){$3.5$}
\put(6.5,1.2){\Large{$S_{13}$}}
\put(11,0.75){\line(-2,-1){1.5}}
\put(11,0.75){\circle*{0.1}}
\put(9.5,0){\circle*{0.1}}
\put(10.95,0.4){$s_{23}$}
\put(9.3,-0.3){$s_{1}$}
\put(10.2,0.5){$3$}
\put(9.5,1.2){\Large{$S_{23}$}}
\end{picture}
\caption{The set $S$ and the three `merged' sets $S_{12}$, $S_{13}$ and $S_{23}$.} \label{c99}
\end{figure}
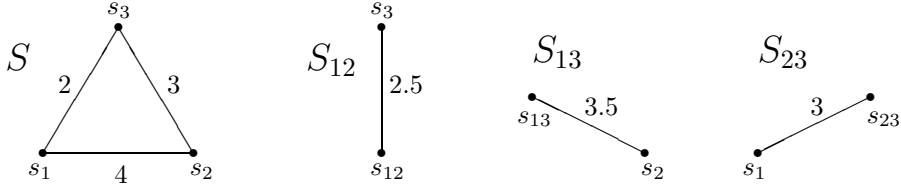

For example,
consider the set $S = \{s_{1},s_{2},s_{3}\}$ illustrated in Figure~\ref{c99},
for which the sets $S_{12}$, $S_{13}$ and $S_{23}$ are also depicted.
Here, we would have $p_{12} = \frac{\frac{1}{4}}{\frac{1}{2} + \frac{1}{3} + \frac{1}{4}} = \frac{3}{13}$
and, similarly,
$p_{13}= \frac{6}{13}$ and $p_{23} = \frac{4}{13}$.
Hence, we would obtain
$D(S) = \frac{3}{13} \left(4 + \frac{5}{2}\right) + \frac{6}{13} \left(2 + \frac{7}{2}\right) + \frac{4}{13} (3+3) = \frac{153}{26}$.

For sets of size three,
this method simplifies to a nice formula:

\begin{Theorem} \label{3setthm}
The diversity measure $D$ defined in equation~(\ref{nice diversity}) satisfies the formula
\begin{displaymath}
D(\{s_{1},s_{2},s_{3}\}) = \frac{3}{2} \left( \sum_{1 \leq i < j \leq 3} \frac{1}{D(\{s_{i},s_{j}\})} \right)^{-1}
+ \frac{1}{2} \sum_{1 \leq i<j \leq 3} D(\{s_{i},s_{j}\}).
\end{displaymath}
\end{Theorem}
\textbf{Proof}
By equation~(\ref{nice diversity}), we have
\begin{eqnarray*}
D(\{s_{1},s_{2},s_{3}\}) &  = & p_{12} \left( D(\{s_{1},s_{2}\}) + \frac{ D(\{s_{1},s_{3}\}) +  D(\{s_{2},s_{3}\}) }{2} \right) \\
& & + p_{13} \left( D(\{s_{1},s_{3}\}) + \frac{ D(\{s_{1},s_{2}\}) +  D(\{s_{2},s_{3}\}) }{2} \right) \\
& & + p_{23} \left( D(\{s_{2},s_{3}\}) + \frac{ D(\{s_{1},s_{2}\}) +  D(\{s_{1},s_{3}\}) }{2} \right) \\
& = & \frac{p_{12}D(\{s_{1},s_{2}\})}{2} + \frac{p_{12}}{2} \sum_{1 \leq i<j \leq 3} D(\{s_{i},s_{j}\}) \\
& & + \frac{p_{13}D(\{s_{1},s_{2}\})}{2} + \frac{p_{13}}{2} \sum_{1 \leq i<j \leq 3} D(\{s_{i},s_{j}\}) \\
& & + \frac{p_{23}D(\{s_{1},s_{2}\})}{2} + \frac{p_{23}}{2} \sum_{1 \leq i<j \leq 3} D(\{s_{i},s_{j}\}) \\
& = & \frac{p_{12}D(\{s_{1},s_{2}\})}{2} 
+ \frac{p_{13}D(\{s_{1},s_{2}\})}{2} 
+ \frac{p_{23}D(\{s_{1},s_{2}\})}{2} \\
& & + \frac{1}{2} \sum_{1 \leq i<j \leq 3} D(\{s_{i},s_{j}\}), \\
& & \textrm{ since } p_{12} + p_{13} + p_{23} = 1 \\
& = & \frac{3}{2} \left( \sum_{1 \leq i < j \leq 3} \frac{1}{D(\{s_{i},s_{j}\})} \right)^{-1} 
+ \frac{1}{2} \sum_{1 \leq i<j \leq 3} D(\{s_{i},s_{j}\}), \\
& & \textrm{ by definition of } p_{kl}.
~
\phantom{qwerty}
\setlength{\unitlength}{.25cm}
\begin{picture}(1,1)
\put(0,0){\line(1,0){1}}
\put(0,0){\line(0,1){1}}
\put(1,1){\line(-1,0){1}}
\put(1,1){\line(0,-1){1}}
\end{picture}
\end{eqnarray*} 

Note that this is the same expression as that given in equation~(\ref{size3formula}),
and so it follows that this measure certainly works satisfactorily for sets of size three.

The equations produced for sets of size greater than three are much more complicated,
and hence more difficult to analyse,
and unfortunately
it seems difficult to find a way to prove that the measure always satisfies Axioms $1$ and $2$
for larger sets.
All experimental results have been positive, however,
and so it is very much hoped that other mathematicians will explore this measure further.

\section{Concluding remarks}

Although many of the properties required of a diversity measure seem simple,
we have seen that it is not easy to produce one.
In particular, the requirement for such a measure to take into account complicated equidistance/symmetry considerations
(and the fact that it it is not obvious precisely what these should be)
seems to make the problem rather difficult.
Nevertheless, it is hoped that the ideas presented in this paper have helped
towards building a rigorous framework for diversity,
and developing a measure that is truly satisfactory.

\end{document}